\documentclass[preprint,letterpaper]{elsarticle}

\usepackage[margin=1.0in]{geometry}
\usepackage{color}
\usepackage{hyperref}
\usepackage{amsmath}
\usepackage{cleveref}
\usepackage{amsfonts}
\usepackage{graphicx}
\usepackage{float}
\usepackage{caption}
\usepackage{subcaption}
\usepackage{fancyhdr}
\usepackage{wrapfig}
\usepackage{booktabs}

\biboptions{square,comma,sort&compress}

\newcommand{\comp}{\mathbb{C}}
\newcommand{\VEC}{\text{v}}
\newcommand{\MAT}{\text{m}}
\newcommand{\Tr}{\text{Tr}}
\newcommand{\Order}{\mathcal{O}}

\begin{document}

\title{Automated Synthesis of Quantum Algorithms\\via Classical Numerical Techniques}
\author[cs]{Yuxin Huang}
\ead{yuxin.huang@vanderbilt.edu}
\author[me]{Benjamin E. Grossman-Ponemon}
\ead{bgrossmanponemon@jcu.edu}
\author[cs]{David A. B. Hyde\corref{cor1}}
\cortext[cor1]{Corresponding author}
\ead{david.hyde.1@vanderbilt.edu}

\address[cs]{Department of Computer Science, Vanderbilt University, 1400 18\textsuperscript{th} Avenue South, Nashville, TN 37212-2846, USA}
\address[me]{Department of Physics, John Carroll University, 1 John Carroll Boulevard, University Heights, OH 44118-4538, USA}

\begin{abstract}

We apply numerical optimization and linear algebra algorithms for classical computers to the problem of automatically synthesizing algorithms for quantum computers.
Using our framework, we apply several common techniques from these classical domains and numerically examine their suitability for and performance on this problem.
Our methods are evaluated on single-qubit systems as well as on larger systems.
While the first part of our proposed method outputs a single unitary matrix representing the composite effects of a quantum circuit or algorithm, we use existing tools---and assess the performance of these---to factor such a matrix into a product of elementary quantum gates.
This enables our pipeline to be truly end-to-end: starting from desired input/output examples, our code ultimately results in a quantum circuit diagram.  We release our code to the research community \textit{(upon acceptance)}.

\end{abstract}


\maketitle

\section{Introduction}

\begin{figure}[H]
\centering
\raisebox{-0.5\height}{\includegraphics[width=0.25\textwidth]{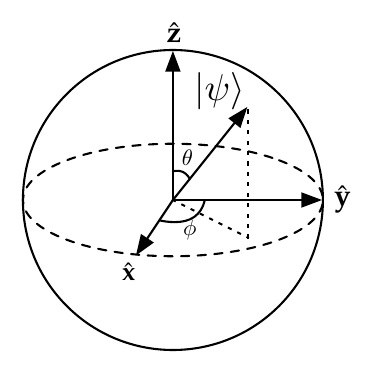}}
\hspace*{.2cm}
\raisebox{-0.5\height}{\includegraphics[width=0.45\textwidth]{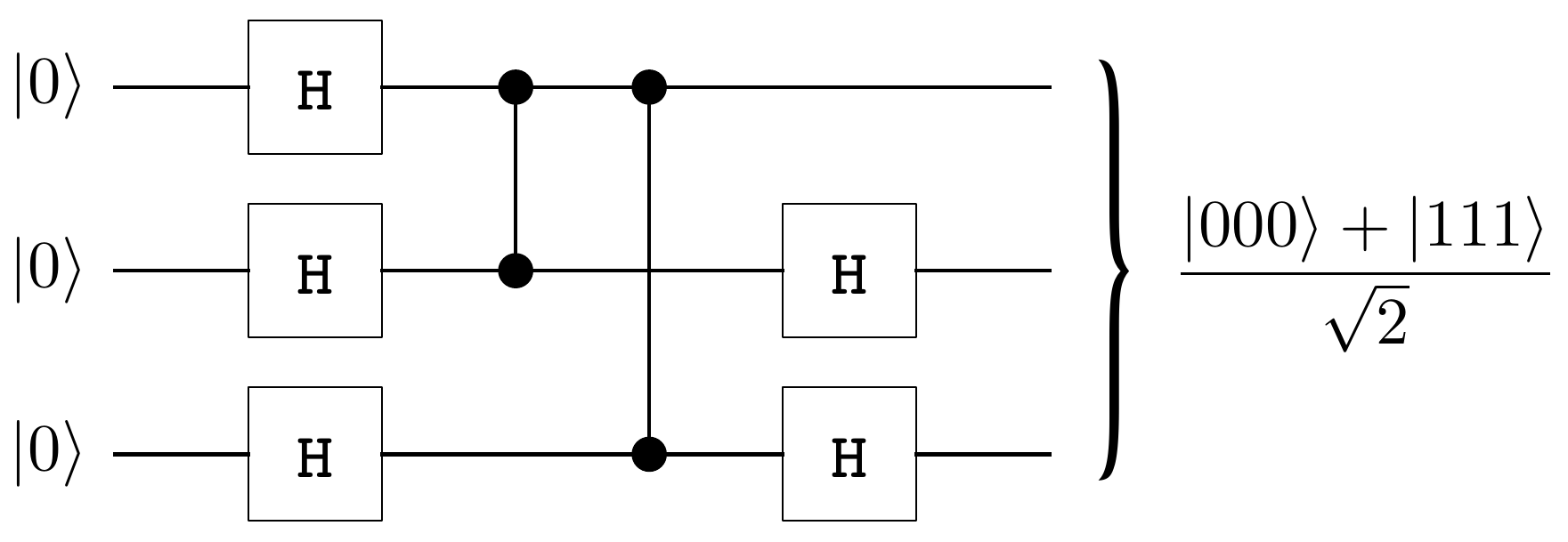}}
\caption{{\em Left}: the Bloch sphere.  $|\Psi\rangle$, the state of a single qubit, can be any point on the surface of the sphere, meaning any superposition of the $|0\rangle$ (+$\hat{\textbf{z}}$) and $|1\rangle$ (-$\hat{\textbf{z}}$) states.
{\em Right}: a simple quantum circuit diagram.  In this algorithm, Hadamard transforms and controlled phase gates are applied to a three-qubit system to create the simplest Greenberger-Horne-Zeilinger state \protect\citep{greenberger_going_1989}.}
\label{fig:quantum_stuff}
\end{figure}

A quantum computer is a system based on {\em qubits}, discrete computational particles that can be in any superposition of ``$0$'' and ``$1$'' states (denoted $|0\rangle$ and $|1\rangle$) at any given
time, unlike classical bits, which must be wholly either $0$ or $1$.  The possible states of a qubit may be visualized
as points on the surface of the Bloch sphere, shown in \Cref{fig:quantum_stuff}.  Based on this physical property, a $k$-qubit system can represent some probabilistic combination of $2^k$ states at
once, an
exponential memory advantage over classical computers.

The development of algorithms for quantum computers dates back to the 1980s (see \citep{nielsen_quantum_2000} for a brief historical discussion).  Since that time, several famous quantum algorithms have
been invented that demonstrate the theoretical power of quantum computers; for example, Shor's algorithm for factoring integers and Grover's algorithm for search \citep{shor_polynomial-time_1997,nielsen_quantum_2000}.  While many successful quantum algorithms have
been proposed by theoreticians, the discovery of new quantum algorithms for increasingly challenging theoretical problems remains accordingly difficult.  Classical numerical techniques may provide one
way to ease this challenge, through the automated discovery, or learning, of quantum algorithms.

In this work, we present a framework for using classical optimization and linear algebra methods to automatically synthesize quantum algorithms merely from sample input and output states of a system. 
Such sample states may be any valid states of a quantum computer, so a theoretician may specify desired input-output pairs, and our framework returns a quantum algorithm that meets the desired behavior (as closely as possible, as we will discuss).
In principle, our work generalizes to arbitrarily large quantum computers, as well as to a more general mathematical context of learning unitary operators between sets of vectors.
Notably, although the first stage of our pipeline outputs a single unitary matrix, we leverage recent advances in the literature to factor the matrix into a product of elementary quantum gate matrices.  Thus, our pipeline allows users to obtain an accurate quantum circuit diagram merely from specifying desired inputs and outputs of an algorithm, which to the authors' knowledge is a novel contribution.

\section{Background and Related Work}

The mathematical formalism of quantum mechanics motivates the problem we consider in this paper.  The state of a $k$-qubit system may be represented by a vector $|\Psi\rangle \in \comp^{2^k}$, where each entry of $|\Psi\rangle$ corresponds to
the probability amplitude of $|\Psi\rangle$ being in the corresponding one of $2^k$ possible states of the system (the ``basis states'').  Just as states may be represented by complex vectors, operations on quantum computers may be represented by complex matrices, which form the most natural set of endomorphisms on $\comp^{2^k}$.  In 
particular, quantum mechanics postulates that quantum systems evolve in unitary ways \cite{qmc062}; that is, operations on quantum computers may be described by {\em unitary} matrices called quantum gates.
Analogously 
to classical computers and circuits, quantum algorithms may be implemented using a sequence (mathematically, a product) of quantum gates.  We therefore restrict our attention in this paper to unitary matrices.

Since the product of unitary matrices is unitary, the problem of synthesizing a quantum algorithm---i.e., a product of gates---can then be viewed as learning a unitary matrix that maps a set of training input vectors to a set of training output vectors with minimum error.
Suppose for a $k$-qubit system we have $m$ given input states $X = [x^{(1)}|\cdots|x^{(m)}] \in \comp^{2^k \times m}$ and corresponding output states $Y = [y^{(1)}|\cdots|y^{(m)}] \in \comp^{2^k
\times m}$. We seek a unitary transformation $U : \comp^{2^k} \rightarrow \comp^{2^k}$ such that $UX \approx Y$. In particular, we seek the operator $U$ that minimizes the Frobenius norm error,
\begin{equation}\label{eq:OrthogonalProcrustesObjective}
f(U) = \frac{1}{2} \| UX - Y \|_F^2.
\end{equation}
The problem of minimizing $f(U)$, if $U$, $X$, and $Y$ are real, is the well-known Orthogonal Procrustes problem \citep{schonemann_generalized_1966}, and the analytical solution is given by $U = WZ^H$, where $W$ and $Z$ arise
from the singular value decomposition $YX^H = W \Sigma Z^H$ (here $H$ denotes the Hermitian, or conjugate transpose, of a matrix).  While this analytical solution solves this basic variation of the problem, we are interested in more generalized versions of this problem
where variational techniques should or must be used.  For example, one may be interested in arbitrary norms, not just the Frobenius norm.  Moreover, with additional constraints or noise in the data (e.g., $U$, $X$, and $Y$ may all be independently subject to noise), it may be the case that the $U$ found through optimization may not be unitary. Thus, we introduce the unitarization error,
\begin{equation}\label{eq:UnitarizationObjective}
g(U) = \frac{1}{4} \| U^H U - I \|_F^2,
\end{equation}
where $U^H$ denotes the Hermitian (or conjugate transpose) of the matrix $U$. Equation \ref{eq:UnitarizationObjective} may be seen as a constraint on the operator $U$ or as a second objective.  We may
further generalize the problem we consider by adding additional constraints, such as placing restrictions on the matrix's condition number, sparsity, etc.  Finally, we note that the Orthogonal Procrustes problem requires the computation (or close approximation) of an SVD, which is a prohibitively expensive calculation ($\Order(n^3)$) for a large system.
In summary, using more general optimization-based techniques allows us to place arbitrary constraints on our problem and to consider much larger physical scenarios than would be solvable with a simple Procrustes problem.

When running optimization-based algorithms on this class of problems, an important property is the number of training examples $m$. If $m<2^k$, then there exists an infinite number of
possibilities for the unitary transformation $U$---the problem is underconstrained.  Meanwhile, $m > 2^k$ overconstrains the problem.
Interestingly, our numerical results suggest that some degree of underconstraining or overconstraining may positively affect convergence.

To the authors' knowledge, attempts to learn quantum algorithms with classical optimization techniques were first made by Ventura \citep{ventura_learning_2000}.
In that work, the author considered a rudimentary delta learning rule for the problem at hand and presented a single example of his algorithm.
In the follow-up work \citep{toronto_learning_2006},
a gradient descent approach was considered for the same problem, with a handful of examples.
In the present work, we compare the costs and convergence rates of previously proposed and newly considered algorithms as we vary the problem size.
In addition, we investigate alternate constraints, additional constraints, and various means of satisfying them.
Per the previous paragraph, we also investigate the impact on our algorithms of varying the number of training examples.
We demonstrate our methods on more complex examples and provide high-performance implementations.
We note that, while the present work is most closely inspired by \citet{ventura_learning_2000} and \citet{toronto_learning_2006}, other authors have performed similar studies using techniques like genetic algorithms \citep{giraldi_ga,giraldi2004genetic}, particle swarm optimization \citep{vlachogiannis2010learning,vlachogiannis2013improved}, and quantum adiabatic optimization \citep{liu2014learning}.

Several other works have considered similar problems.
For instance, \citet{mahmud2021optimizing} investigated synthesizing initial quantum states from classical data (i.e., learning a vector, not a matrix).
\citet{song2003computational} also studied state synthesis, from an initial zero state, via unitary factorization in time $O\left(2^{2k}\right)$.
\citet{10.5555/2011670.2011675} provides a method for transforming an arbitrary quantum state $|a\rangle$ to another arbitrary state $|b\rangle$, but their method is designed for just one input-output pair, and they compute a specific structure of quantum circuit that may not be optimal.

\section{Estimating a Single Unitary}

The first part of our proposed pipeline is to determine a single unitary matrix that maps given input data to given output data as accurately as possible; we will later discuss factoring this matrix into a product of quantum gates.
For a $k$-qubit system, we seek unitary operators $U \in \comp^{2^k \times 2^k}$, meaning that we have to find $2^{2k}$ unknown complex values. For this reason, we primarily investigated gradient-based learning algorithms.
Hessian-based algorithms like Newton's method may require fewer total iterations, but the computation of the Hessian (in general, a dense array of $2^{4k}$ complex entries) and its subsequent inversion is likely to be prohibitive for larger systems.

Within gradient-based algorithms, we considered batch gradient descent and a single-layer neural network updated via the delta learning rule (effectively equivalent to stochastic gradient descent).
Incorporating constraints like Equation \ref{eq:UnitarizationObjective}, we investigated the impact of penalty methods, Lagrange multiplier techniques, and variable learning rates on these gradient methods.  

For the fully constrained (or overconstrained) noiseless problem (i.e., $m\geq2^k$ and $Y = UX$), the operator $U$ in the unconstrained setting may be found via minimization of Equation \ref{eq:OrthogonalProcrustesObjective}. The gradient descent rule is given by
\begin{equation}\label{eq:BatchGradientDescent}
U := U - \alpha (UX - Y)X^H ,
\end{equation}
where $(UX-Y)X^H$ is the gradient of Equation \ref{eq:OrthogonalProcrustesObjective}, as derived in \citet{toronto_learning_2006}, and $\alpha$ is the learning rate. 

Incorporating the unitarization constraint, Equation \ref{eq:UnitarizationObjective}, for under-constrained problems (or problems where one wishes to ensure $U$ is unitary) may be handled through penalties or Lagrange multipliers. Toronto and Ventura use a penalty formulation of the problem with an additional unitarization step \citep{toronto_learning_2006}. During the first phase, the update rule is
\begin{equation}\label{eq:BatchGradientDescentPenalty}
U := U - \alpha (UX - Y)X^H - \beta U(U^HU - I) ,
\end{equation}
where $U(U^H U-I)$ is the gradient of Equation \ref{eq:UnitarizationObjective}. Incorporating the unitarization constraint with Lagrange multipliers gives rise to the update rules
\begin{equation}\label{eq:BatchGradientDescentLagrange}
\begin{aligned}
U &:= U - \alpha (UX - Y)X^H - \lambda U (U^HU-I), \\
\lambda &:= \lambda + \frac{\beta}{4} \|U^H U - I\|_F^2.
\end{aligned}
\end{equation}

The above algorithms rely on handling all the training data at once. As an alternative, stochastic gradient descent (or the delta learning rule \cite{ventura_learning_2000}) may be used. The update
rules are identical in form as those above, but only operate on a single training example at a time. For example, without unitarization, the update rule is given by
\begin{equation}\label{eq:StochasticGradientDescent}
U := U - \alpha(Ux^{(i)} - y^{(i)})(x^{(i)})^H.
\end{equation}

Despite its likely limited scalability, Newton's method was also investigated.
The difficulty in implementing this method lies in that the Hessian of the objective, Equation \ref{eq:OrthogonalProcrustesObjective}, is a fourth-order tensor.
To address this, let $\VEC(A)$ be a re-indexing of a matrix $A$ into a vector, and $\MAT(B)$ be the re-indexing of a fourth-order tensor $B$ into a matrix.
Then the update rule for Newton's method is
\begin{equation}\label{eq:NewtonsMethod}
\VEC(U) := \VEC(U) - \alpha\, \MAT(\nabla_U^2 f(U))^{-1} \VEC((UX-Y)X^H).
\end{equation}
The form of the Hessian of Equation \ref{eq:OrthogonalProcrustesObjective}, in the general complex, rectangular setting, is derived in \ref{sec:Hessian}; the authors have not found this mathematical derivation in prior literature.

\subsection{Sequential Matrix Learning}
\label{sec:sequential}

Each entry in a vector $y^{(i)} = U x^{(i)}$ is independently determined by a different row of $U$.
Motivated by this, we propose a novel (to the authors' knowledge) formulation of the target optimization problem in order to achieve a more efficient algorithm.
Equation \ref{eq:OrthogonalProcrustesObjective} can be split into a sum of norms over the rows of the matrices $U$ and $Y$,
$$
f(U) = \sum \limits_{j=1}^{2^k} \frac{1}{2} \|e_j^T U X - e_j^T Y\|_F^2,
$$
where $e_j^T U$ is the $j$\textsuperscript{th} row of the matrix $U$. Furthermore, both the batch and stochastic gradient descent update rules can be split into updates for each row.
For the $j$\textsuperscript{th} row, we have
$$
e_j^T U := e_j^T U - \alpha (e_j^T U X - e_j^T Y) X^H.
$$
Thus, we may find the optimal matrix $U$ row-by-row.
This \emph{sequential} formulation is easily parallelized.
For example, on a shared-memory architecture, we may provide each thread with a given row of $U$ and $Y$, and all threads may access the data $X$.
This level of parallelization is on top of any parallelization of BLAS operations like matrix-vector multiplications that may be implemented in a particular linear algebra library.

\subsection{Data Generation}
\label{sec:data-generation}

For testing our methods, it is possible to generate arbitrarily many examples for our framework, rather than limiting ourselves to one or two examples as in prior works \citep{ventura_learning_2000,toronto_learning_2006}.
Briefly, any quantum operator may be approximated arbitrarily well by (possibly very long) products of matrices from a set of {\em universal} gates, e.g., the three gates $H$, $R(\pi/4)$, and $CNOT$ for a two-qubit system.
Based on this idea, one can propose a simple algorithm for generating random quantum algorithms.
Pick a number of qubits $k$ and a random $n\geq1\in\mathbb{Z}$.
For $i=1,\dots,n$, pick a random operator from the universal gate set, and apply it to a random choice of qubit(s) in the system (so, part of the matrix may be an identity block).
The product of the matrices selected at each step generates a random quantum algorithm, which we can then attempt to learn using our framework.
Knowing the algorithm's matrix form allows us to generate exactly compatible training input and output vectors; we can then pretend to ``forget'' the matrix and use these training pairs to learn it, while verifying the result against the true matrix.

However, this idea raises a few questions, such as what constitutes an appropriate choice of $n$.
Recognizing that quantum algorithms can be represented using unitary operators, we instead preferred for our experiments an algorithm that generates random unitary operators directly, without requiring knowledge of quantum gates.
If we sample $\left(N_j\right)_{kl} \sim \mathcal{N}(0,1)$ for $j=1,2$ and
compute:
\begin{align*}
U&=N_1 + iN_2\\
U&=\text{normalizedCols}(U)\\
Q,R&=\text{qr}(U)\\
U&=Q\text{diag}\left(\text{diag}(R)/\text{abs}(\text{diag}(R))\right)
\end{align*}
where $\text{normalizedCols}$ normalizes each column of $U$ (and where division and absolute value are elementwise operations), then we obtain a unitary $U$ that is distributed uniformly according to
the Haar measure \citep{alan_edelman_random_2005,ozols_how_2009}.
As with the other proposed generation technique, we may use this algorithm to generate a random $U$ and $X$, compute $Y=UX$, and then forget $U$ and
attempt to learn it again, in order to have random tests for our numerical techniques.
In our tests, when generating a random $X$, we sample $\left(\hat N_j\right)_{kl} \sim \mathcal{N}(0,1)$ for $j=1,2$ and set $X = \text{normalizedCols}(\hat N_1 + i \hat N_2)$.

We note that for our randomly-generated data, we observed degraded performance of all the algorithms considered when the condition number of randomly-generated
$X$ was large.  Specifically, convergence was greatly slowed when a near-zero singular value was present in $X$.  To combat this ill effect, we added a check for each random $X$,
regenerating any matrix that had a condition number above $100$ (unless otherwise specified).

\subsection{Numerical Experiments}

We present the results of learning a single unitary using various numerical techniques.
Throughout, ``GDP'' refers to the batch gradient descent penalty method, ``GDLM'' refers to the batch gradient descent Lagrange multiplier method, ``DLR'' refers to the delta learning rule-updated single-layer neural network method, and ``NM'' refers to Newton's method. 
A ``*'' next to a name indicates that the method is constrained with Equation \ref{eq:UnitarizationObjective}.
Unless otherwise specified, a coefficient of $\beta = 0.1$ was used for the unitarization objective.

Furthermore, unless otherwise stated, constant learning rates (step sizes) were used.
We investigated linesearch methods for learning rate selection, e.g., seeking a learning rate $\alpha$ that satisfies the Wolfe conditions:
{\small
\begin{align*}
f(U^{(n)} + \alpha^{(n)} dU^{(n)}) \leq f(U^{(n)}) +c_1 \alpha^{(n)} (dU^{(n)})^H \nabla_U f(U^{(n)}),\\
(dU^{(n)})^H \nabla_U f(U^{(n)} + \alpha^{(n)} dU^{(n)}) \geq c_2 (dU^{(n)})^H \nabla_U f(U^{(n)})
\end{align*}
}
(sufficient descent and sufficient curvature, respectively).
We also considered backtracking linesearch, which repeatedly divides an initial guess for $\alpha$ until the iterate satisfies the
sufficient descent condition.
We refer the reader to \citet{wright_numerical_2006}.
Although these techniques can be used to ensure convergence, we were generally able to use larger constant learning rates that still yielded convergent results in practice.

For the implementation of our experiments, we developed an original code in C++14 (released on publication).
Matrix algebra was performed through the use of the Armadillo library \citep{sanderson_armadillo:_2010,sanderson2019practical}, version 14.0.1.
BLAS routines were provided via OpenBLAS v0.3.21+ds-4.
We used our own implementations of all optimization algorithms, and data generation was implemented as described in the prior subsection. 
Linear algebra routines were parallelized with the built-in support of Armadillo.
Additionally, the sequential formulation of our problem (learning one row at a time) was trivially parallelized with OpenMP.
For our experiments, we ran our code on a workstation with dual AMD EPYC 75F3 CPUs, for a total of up to 64 cores, on a machine running Debian 12.
Our system had 512GB of RAM, though our code used only a few MB even for our largest tests.

Convergence was evaluated using a tolerance of 1e-15.
For each method, we check convergence by evaluating whether the objective value is below the tolerance.
Sometimes convergence stalls; we consider a test to have stalled if the difference in the objective value between the prior and current iteration is less than 1e-18.
Instances where a test stalls are counted separately and do not count towards measurements such as the average number of iterations for a test.
We select sufficiently small step sizes (learning rates) that our results do not include cases where the objective diverges.

\subsubsection{Examples from Literature}

\begin{figure}
    \centering
    \includegraphics[width=0.49\textwidth]{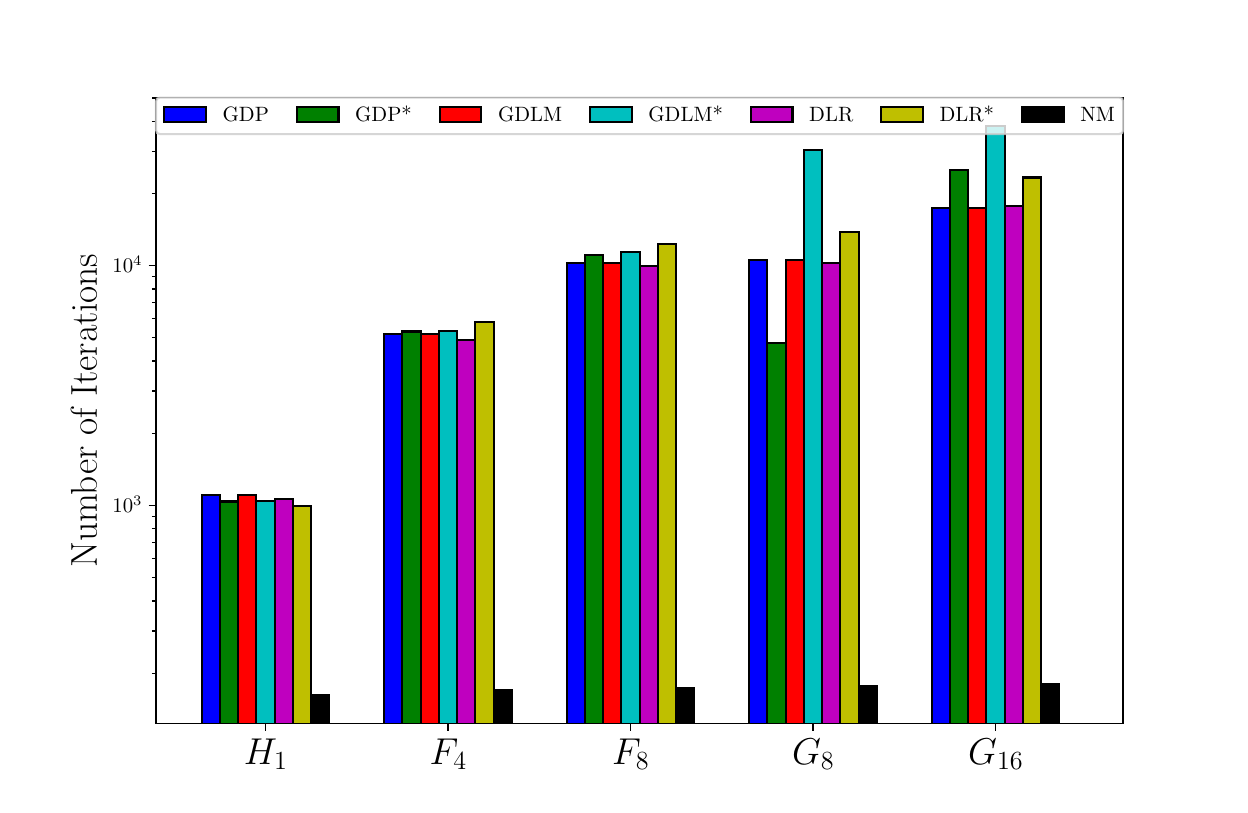}
    \caption{Number of iterations to convergence for sample quantum operators. Since some methods stalled more frequently than others, error bars are scaled differently, so we omit them for fairness. All methods run with constant learning rate $\alpha=0.1$ in an attempt to isolate only differences in learning algorithm.}
    \label{fig:barchart1}
\end{figure}

We tested each of the numerical approaches on the following well-known gates: the Hadamard gate on a single qubit,
\begin{equation}\label{eq:hadamard}
H_1 = \frac{1}{\sqrt{2}} \begin{bmatrix}
1 & 1 \\
-1 & 1
\end{bmatrix},
\end{equation}
the Quantum Fourier Transform of size $N$,
\begin{equation}\label{eq:qft}
(F_N)_{jk} = \left(e^{2\pi i/N}\right)^{(j-1)(k-1)}/\sqrt{N},
\end{equation}
and Grover's iterate of size $N$,
\begin{equation}\label{eq:grover}
(G_N)_{jk} = 2/N - \delta_{jk},
\end{equation}
$j,k=1,\dots,N$.
Shown in Figure \ref{fig:barchart1} are the number of iterations required for each method to converge on random data $X$ and with constant learning rate $\alpha = 0.1$.
The results in the figure were averaged over up to 100 random $X$, and $N$ is indicated in the subscripts along the horizontal axis.
In cases where stalling occurred, we only computed statistics using the subset of trials that did not stall.
Table \ref{tab:lit-examples} quantifies the data shown in Figure \ref{fig:barchart1} and also lists how many trials stalled for each method.

\begin{table}
    \centering
    \begin{tabular}{cccccccc}\toprule
      Test & GDP & GDP* & GDLM & GDLM* & DLR & DLR* & NM \\ \midrule
      $H_1$ & 1109, 2 & 1039, 2 & 1110, 2 & 1040, 2 & 1061, 2 & \textbf{992}, 2 & 162, 0 \\
      $F_4$ & 5167, 5 & 5316, 5 & 5168, 5 & 5363, 5 & \textbf{4917}, 5 & 5824, 5 & 169, 0 \\
      $F_8$ & 10251, 19 & 11116, 19 & 10252, 19 & 11369, 19 & \textbf{9979}, 18 & 12259, 18 & 174, 0 \\
      $G_8$ & 10518, 19 & \textbf{4779}, 22 & 10519, 19 & 30225, 16 & 10239, 18 & 13862, 18 & 177, 0 \\
      $G_{16}$ & \textbf{17431}, 70 & 25001, 76 & 17432, 70 & 38213, 70 & 17729, 67 & 23317, 68 & 181, 0 \\ \bottomrule
    \end{tabular}
    \caption{Statistics for each of the tests shown in Figure \ref{fig:barchart1}, rounded to the nearest integer. For each method, the results reported are pairs of (average number of iterations, number of trials that stalled).  The iteration counts are averaged only over the trials that did not stall.  100 trials are run for each scenario.  The best first-order method result is highlighted in each row.}
    \label{tab:lit-examples}
\end{table}

We include unconstrained Newton's method in these results for illustration, highlighting the benefits of quadratic convergence.
However, as mentioned, Newton's method requires the computation and inversion of the Hessian of Equation \ref{eq:OrthogonalProcrustesObjective}, the size of which grows like $\Order(2^{4k})$ for a $k$-qubit system.
Furthermore, because the Hessian, Equation \ref{eq:OrthogonalProcrustesHessian}, is not very sparse, Newton's method will not be feasible for large system sizes.
Accordingly, we focus on assessing the relative performance of gradient-based methods; the best-performing gradient-based method for each test problem is highlighted in bold in Table \ref{tab:lit-examples}.

The results in Figure \ref{fig:barchart1} suggest that several of the gradient-based methods generally exhibit similar performance.
Notably, we observe that methods constrained with Equation \ref{eq:UnitarizationObjective} often (though not always) required more iterations than their unconstrained counterparts.
We also note that GDLM* appears to perform consistently worse on the two Grover's iterate test cases, which are the most complex of these elementary tests.
We emphasize that by using smaller and/or adaptive learning rates, it would be straightforward to improve the convergence of these methods; the intent of these results is to attempt to isolate the differences in the optimization algorithms' behavior, rather than to showcase each algorithm running at its best.

\subsubsection{Random Unitary Operators}
\label{sec:wall-clock}

\begin{figure}[!ht]
\centering
\includegraphics[width=0.5\textwidth]{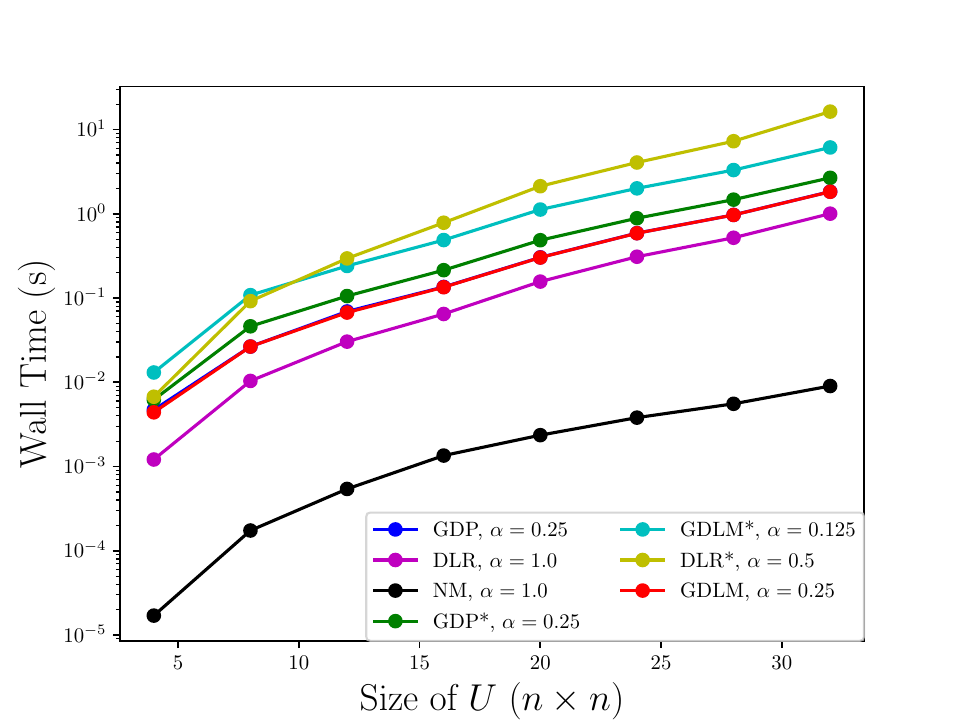}
\caption{Wall time to converge for each method versus size of random system.  Here we aim to show each algorithm at its best---using the largest feasible constant step size---hence some methods appear to significantly outperform others.}
\label{fig:plot1}
\end{figure}

We also tested each algorithm on random unitary operators $U \in \mathbb{C}^{n \times n}$ and with random data $X \in \mathbb{C}^{n \times n}$.
For these tests, we compared the run times against system size, where $n$ ranges from $4$ to $32$.
These matrices were generated according to the second random algorithm presented in Section \ref{sec:data-generation}.
In order to show each method at its best, we manually tuned the learning rate to get peak performance; for each method, we selected a learning rate that yielded convergent results at all problem sizes, starting with $\alpha = 1.0$ and reducing by half until no stalling occurred.
We used constant learning rates for these tests that, in general, exceeded the step sizes suggested by linesearch methods for guaranteed convergence.
Results are shown in Figure \ref{fig:plot1}.

The figure illustrates that certain methods were able to leverage significantly larger learning rates; combined with the observation from Figure \ref{fig:barchart1} that the methods converge in a similar number of iterations, it is thus unsurprising that methods like the delta learning rule and Newton's method have significantly smaller wall clock times.
This suggests that methods such as DLR are more numerically robust, especially since these step sizes are beyond the theoretically-justified step sizes given by linesearch methods.
On the other hand, when constraints are introduced, DLR\textsuperscript{*}---even with a two or four times larger step size---performs comparably to or slower than the other constrained techniques, suggesting that DLR may only have a net advantage for the simpler unconstrained case.

\subsubsection{Number of Training Examples}

\begin{figure}[!ht]
\centering
\includegraphics[width=0.5\textwidth]{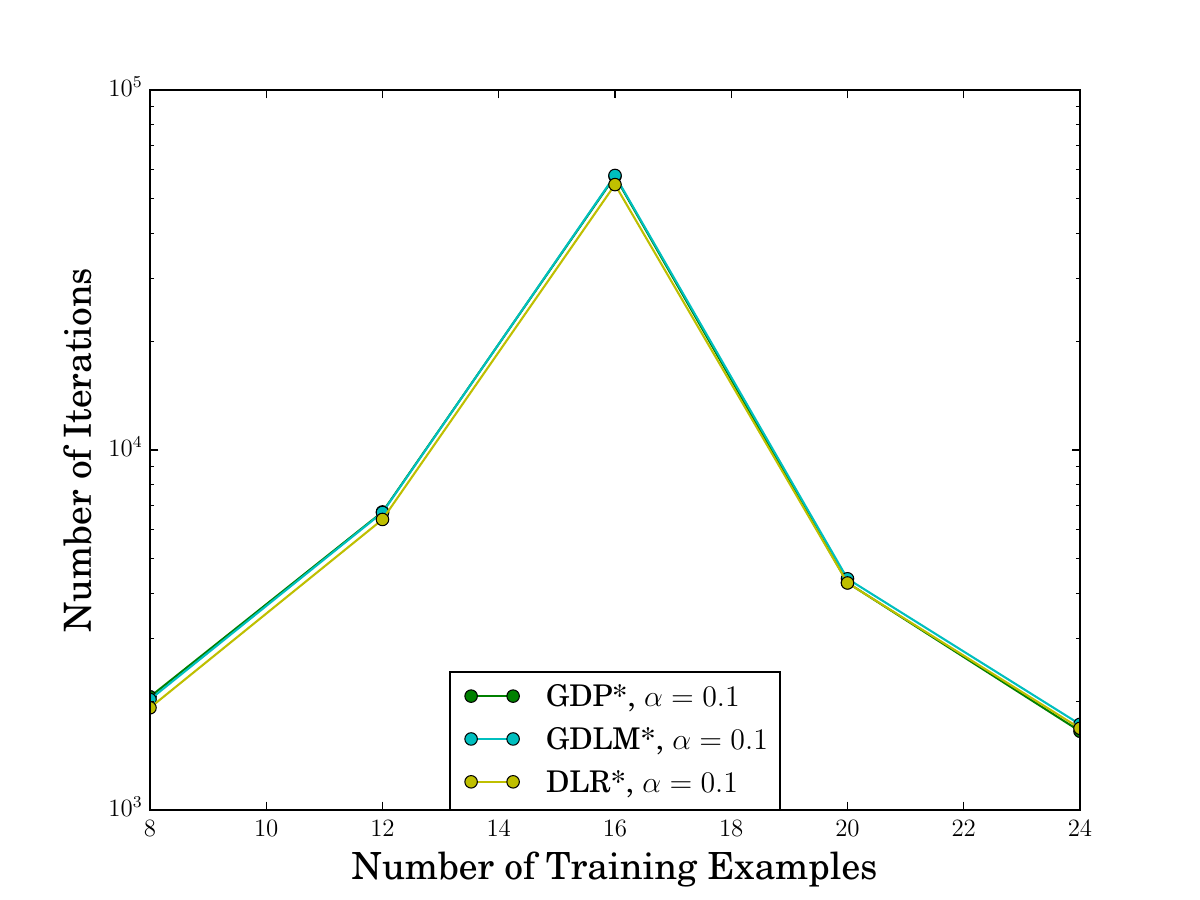}
\caption{Number of iterations to convergence versus number of training examples on a random 4-qubit system.  For this example, underconstraining or overconstraining the problem accelerates convergence.}
\label{fig:plot2}
\end{figure}

In \citet{toronto_learning_2006}, the authors considered the learning of operators for $m < 2^k$ training examples. We sought to compare the performance of the learning algorithms as we varied the
number of training examples, both underconstraining and overconstraining the system. 
Figure \ref{fig:plot2} shows the results of varying the number of examples used to learn a random unitary with $k=4$; we focus on the constrained versions of the algorithms for this figure.
The figure shows that, with underconstrained data, each algorithm was able to converge in fewer steps, likely due to the fact that the solution $U$ is no longer unique (hence increasing the solution space).
Meanwhile, overconstraining also reduced the number of iterations to convergence.
While the additional training examples are no longer linearly independent, the additional examples were consistent, since we generated the example outputs as $Y=UX$ (using the exact $U$).
Thus, intuitively, the additional information appears to help the algorithms reach the optimal solution in fewer iterations. 
However, the addition of more training examples increases the computational cost of each iteration. 
For example, with batch gradient descent, we must compute the update $(UX - Y)X^H$.
If $X \in \comp^{2^k \times m}$, then the number of operations involved in the matrix-matrix products is $\Order(2^{2k}m)$. 
Hence, while additional consistent examples do appear to aid in convergence, larger $m$ have potential to significantly slow down the considered algorithms.

\subsubsection{Sequential versus Batch Formulations}

\begin{figure}[H]
\centering
\includegraphics[width=0.5\textwidth]{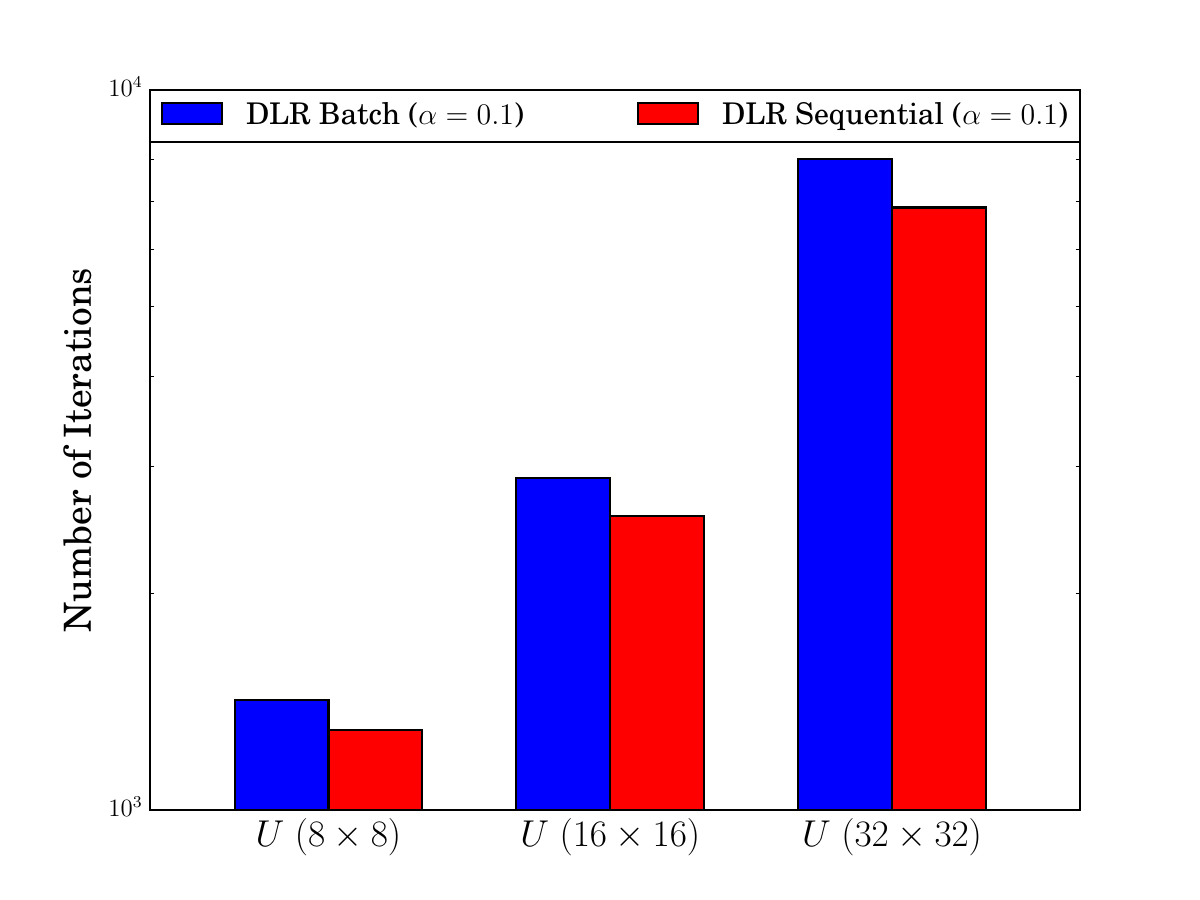}
\caption{Number of iterations to converge for batch versus sequential formulations of gradient descent.  Results for the sequential formulation are normalized by $n$, i.e., we show iterations to learn one row, because we implement finding each row in parallel and because an iteration of the sequential formulation takes about $1/n$ of the time of an iteration of the batch formulation.  DLR was selected for these tests due to its advantageous performance in Section \ref{sec:wall-clock}.}
\label{fig:BarChart2}
\end{figure}

We also compared the performance of our sequential (row-by-row) reformulation of stochastic gradient descent against the batch (or full $U$) version on three random unitary operators $U$ of varying sizes. 
As seen in \Cref{fig:BarChart2}, it takes fewer iterations to learn a row of $U$ in the sequential formulation than it does to learn all of $U$ at once in the batch formulation.
This demonstrates that the sequential formulation is, as the theory in Section \ref{sec:sequential} suggests, solving an easier optimization problem than that solved in the batch formulation.
The number of unknowns in the sequential formulation scales like $\Order(n)$, not $\Order\left(n^2\right)$ as in the batch formulation.
Thus, as $n$ continues to grow, the batch formulation soon becomes infeasible---smaller and smaller learning rates are required, and convergence becomes challenging to obtain.
We believe that our sequential formulation offers a path towards synthesizing  operators for much larger quantum systems.

\subsubsection{Input-Output Distance}

\begin{figure}[!ht]
\centering
\includegraphics[width=0.5\textwidth]{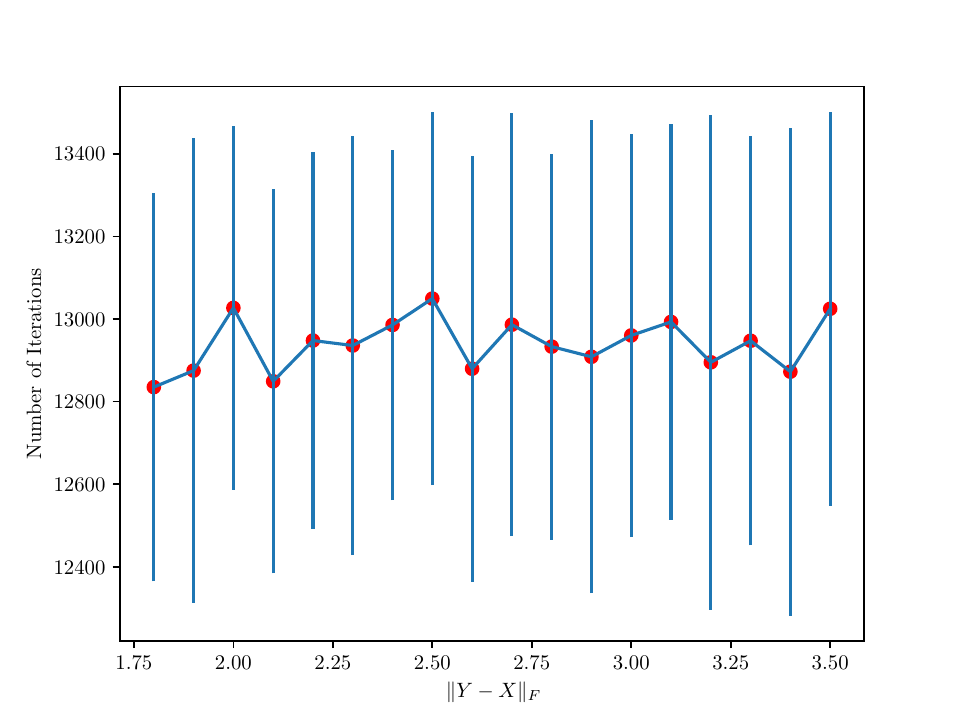}
\caption{The effect of the Frobenius norm distance between $X$ and $Y$ on convergence of the GDP\textsuperscript{*} method.  Mean and standard deviation error bars are shown for each bin of samples (bin width of 0.1).  No significant correlation is observed ($R^2$ coefficient = 0.23).}
\label{fig:xydiff}
\end{figure}

We considered further factors that may impact convergence.
The first factor is the distance between input and output states in the Frobenius norm, $\|Y-X\|_F$.
Intuitively, one may expect that if $Y$ is very close to $X$, learning the mapping operator $U$ requires fewer iterations, since the true $U$ will be closer to the identity matrix and may have lower condition number.
Using the penalty formulation of the gradient descent approach ($\alpha=\beta=0.1$), we solved the constrained sequential formulation of the problem 50,000 times for random square $X$ and $U$ with $n=4$, recording $\|Y-X\|_F$ and the number of iterations needed for convergence.
\Cref{fig:xydiff} shows the result of our experiment.
We note that in order to obtain roughly equal numbers of samples for different values of $\|Y-X\|_F$, we split the range of observed $\|Y-X\|_F$ into bins of width
$0.1$ and only used the first 100 samples in each bin for our analysis.
This ensures the results are unbiased by the differing frequencies of observed $\|Y-X\|_F$ during random trials.
The figure shows that there is no clear trend between the number of iterations required for convergence and $\|Y-X\|_F$.
The $R^2$ coefficient of the best fit line through the bin means is 0.23, indicating that there is no significant dependence on $\|Y-X\|_F$ for convergence.

\subsubsection{Initial Guess Selection}

\begin{figure}[!ht]
\centering
\includegraphics[width=0.5\textwidth]{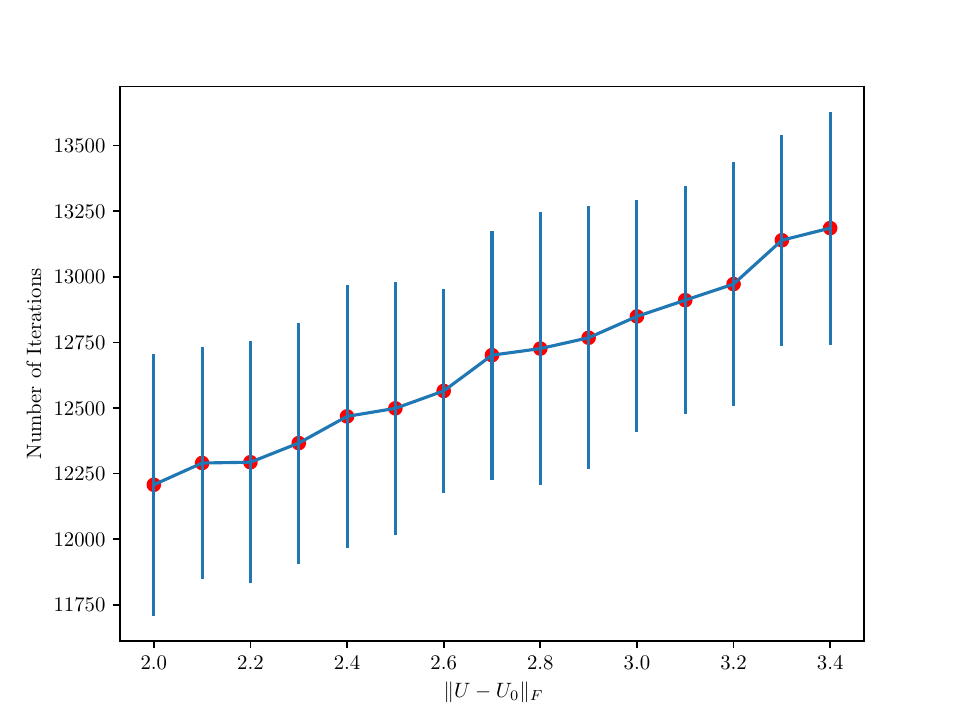}
\caption{The effect of the Frobenius norm distance between $U$ and the initial guess $U_0$ on convergence of the GDP\textsuperscript{*} method.  Mean and standard deviation error bars are shown for each bin of samples (bin width of 0.1).  An $R^2$ coefficient of 0.99 is computed for the line of best fit through the bin means, indicating a statistically significant linear fit.}
\label{fig:uu0diff}
\end{figure}

On the other hand, we observe that algorithms may exhibit a significant linear dependence between the number of iterations and the quality of the initial guess for $U$, which we quantify as $\| U - U_0 \|_F$.
Figure \ref{fig:uu0diff} evaluates the number of iterations the GDP\textsuperscript{*} method (as in the prior subsection) takes to converge as the quality of the initial guess becomes worse.
$n=4$ as before.
As in the prior subsection, we bin results to intervals of width 0.1 and consider the first 100 sample points within each bin.
We find that the line of best through the bin means has an $R^2$ coefficient of 0.99, indicating a strong linear dependence on the quality of the initial guess.

\subsubsection{Conditioning of \texorpdfstring{$X$}{X}}

\begin{figure}[H]
\centering
\includegraphics[width=0.5\textwidth]{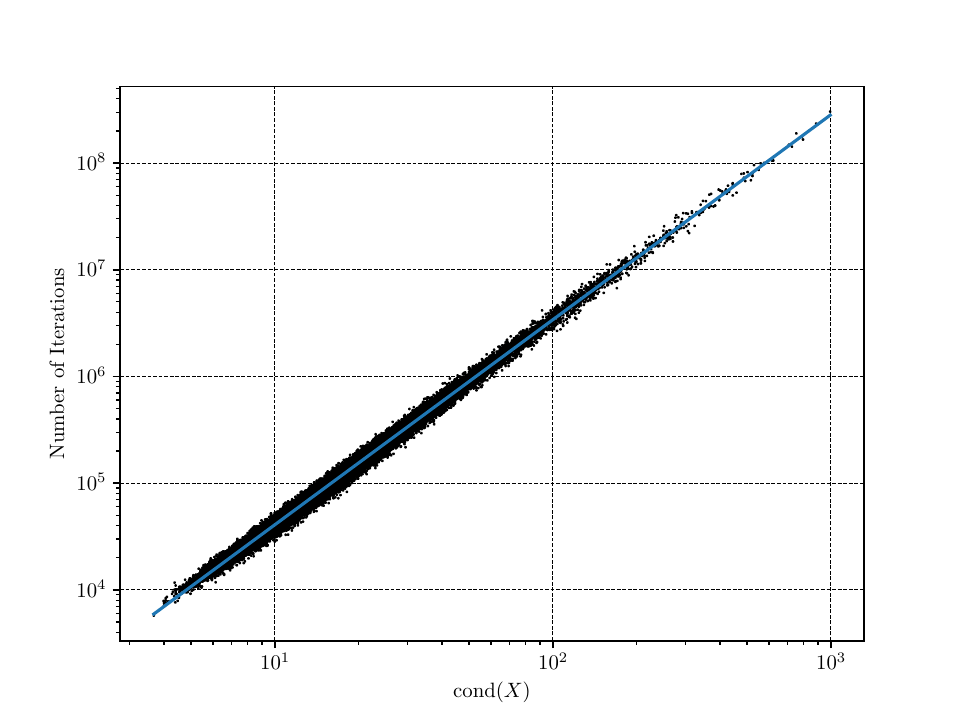}
\caption{The effect of $\text{cond}(X)$ on the convergence of the GDP\textsuperscript{*} method. Plotted with log-log scaling, a statistically significant linear fit (with $R^2$ coefficient of $0.99$ and slope $1.92$) is evident, suggesting that the condition number of $X$ can make optimization algorithms require quadratically more iterations to converge.}
\label{fig:xcond}
\end{figure}

We find that a greater arbiter of convergence speed is the condition number of $X$.  Since unitary matrices have condition number $1$, the conditioning of the problem is largely determined by $X$ (though, of
course, it will likely be affected by added constraints).  We generate 50,000 random $X$ matrices ($n=4$ as before) and use GDP\textsuperscript{*} to measure the number of iterations it takes to converge ($U$ is also randomized for each trial).
Figure \ref{fig:xcond} shows the results.
On a log-log scale, the fit is almost perfectly linear, with a slope of $1.92$ an $R^2$ coefficient of 0.99.
Thus, the condition number of $X$ can quadratically increase the number of iterations required for an optimization method to converge on our target problem.

\section{Factorization of Operators into Gates}

In order to build an end-to-end pipeline for synthesizing quantum algorithms, we now turn our attention to factoring unitaries learned from the techniques of the prior section into products of quantum gates.
Although we leverage existing techniques for this purpose, our evaluation of several recent methods on random unitaries provides helpful insights into the state of the art.

As quantum computing's promise and popularity have increased over the past several decades, increasing attention has been given towards developing unitary factorization algorithms.
Perhaps the most well-known is the Solovay-Kitaev algorithm and associated theorem \citep{nielsen_quantum_2000,dawson2005solovay,kitaev2002classical}, which proves that with enough quantum gates, an arbitrary single-qubit unitary can be approximated arbitrarily well.
Around the same time, work by \citet{PhysRevA.52.3457} and by Tucci \citep{tucci1998rudimentary,tucci1999rudimentary} considered differing approaches to unitary factorization, with the former based on the recursive algorithm of \citet{reck1994experimental} (see also Chapter 2 of \citet{murnaghan1958orthogonal}) and the latter based on the CS decomposition \citep{stewart1982computing,van1985computing,paige1994history}.
Since these early algorithms, research has largely focused on producing methods with improved efficiency (either faster performance or reduced upper bounds on resulting circuit complexity).
For instance, \citet{10.1145/3505181} investigate approximately factoring a unitary matrix into a quantum circuit and demonstrate that the Quantum Shannon Decomposition \cite{shende2005synthesis} can be compressed without practical loss of fidelity.
\citet{rudolph2022decomposition} focus on achieving shallow parameterized quantum circuits when factoring matrix product states.
The OpenQL library \citep{openql2021} library uses the recent factoring algorithm of \citet{krol2022efficient}, which is based on the Quantum Shannon Decomposition.
Another open-source codebase, SQUANDER, leverages the recent algorithms of \citet{rakyta2022approaching}, \citet{rakyta2022efficient}, and \citet{rakyta2022highly}.
We also highlight the work of \citet{GOUBAULTDEBRUGIERE2020107001}, which performs unitary factorization using the QR factorization via Householder transformations, enabling efficient multicore CPU and GPU implementations.
Qiskit, as a comprehensive open-source quantum library, incorporates a variety of decomposition methods within its transpiler.
At its core, Qiskit uses Quantum Shannon Decomposition for unitary synthesis.
After this, Qiskit can perform several optimization passes in order to reduce the complexity of the generated circuit.
At optimization level 1, Qiskit uses inverse gate cancellation \cite{qiskit_inverse_cancellation}, which removes pairs of gates that are inverses of each other, or single-qubit gate optimization \cite{qiskit_1q_optimization}, which potentially combines chains of single-qubit gates into a single gate. 
Optimization level 2 applies commutative cancellation \citep{qiskit_commutative_cancellation} to cancel redundant (self-adjoint) gates through commutation relations, thus reducing the total gate count.
At level 3, Cartan's KAK Decomposition (described in \citet{tucci2005KAKDecomp}) is used for resynthesizing two-qubit blocks, which helps to reduce the number of CNOT gates in the circuit.
However, it is important to note that many of the optimization techniques in Qiskit's transpiler are heuristic-based, meaning that spending more computational effort does not always guarantee an improvement in the quality of the output circuit.

\subsection{Numerical Experiments}

Given the dynamic state of the field, we focused our attention on evaluating several recent techniques for which source code was readily available: 
(1) Qiskit \cite{qiskit2024} with three distinct optimization levels, (2) OpenQL \cite{openql2021}, and (3) Squander \cite{rakyta2022efficient}. 
We used Qiskit v1.2.0, OpenQL v0.12.0, and Squander v1.7.4.
We evaluate these methods, as feasible, on single- and multi-qubit systems.

We note that there are four main metrics used to evaluate the performance of this class of factorization algorithm: wall clock time, number of resulting CNOT gates, total number of resulting gates, and fidelity.
For instance, the method of \citet{Iten_2016} achieves a lower bound on the number of CNOT gates produced, but may not be the most performant in terms of wall clock time.
Considerations aside from wall clock time are important because the ultimate goal of a pipeline like that proposed in this work is to realize real-world quantum circuits, which should be designed using as few gates as possible given that even the most powerful circuit-model quantum computers available today have severely limited qubit counts and coherence times.
Despite hardware constraint considerations, we are also interested in measuring how well our synthesized circuit reproduces the actions of the original unitary operation. 
Specifically, we use process fidelity \cite{Gilchrist_2005, PhysRevResearch.3.033031}, which measures how similar two quantum operations are to each other.
The process fidelity between two quantum operations $\mathcal{E}$ and $\mathcal{F}$ is mathematically defined as:
\[F(\mathcal{E}, \mathcal{F}) = \frac{1}{d^2} \sum_{i=1}^{d^2} \text{Tr}[\mathcal{E}(\sigma_i^\dagger)\mathcal{F}(\sigma_i)] , \]
where $d$ is the dimension of the Hilbert space, and ${\sigma_i}$ is a complete set of mutually orthonormal operators on this space. This formula quantifies the similarity between the two operations across all possible input states.
Fidelity values range from 0 to 1. A fidelity of 1 indicates perfect reproduction of the original unitary operation, while lower values suggest deviations in the decomposed circuit from the original unitary operation.
In our analysis, we report the error as a direct measure of the accuracy lost during the decomposition process, calculated as $\text{error} = 1 - \text{fidelity}$. That is, a lower error (closer to 0) indicates a more accurate decomposition. 

\begin{figure}
    \centering
    \includegraphics[width=1.0\linewidth]{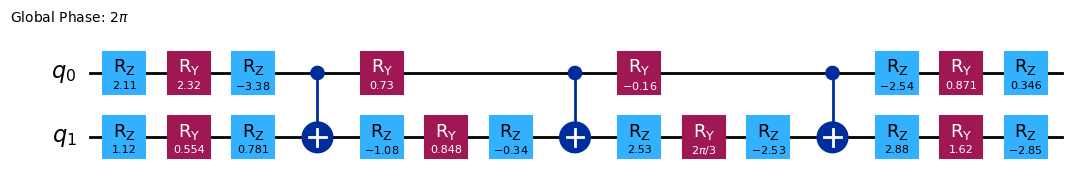}
    \caption{An example of a two-qubit gate decomposition using the basis gate set $R_Z$, $R_Y$ and $CNOT$.}
    \label{fig:decomposition_visualization}
\end{figure}

In order to make a fairer comparison between the algorithms mentioned above, we attempted to use the same universal basis gate set for each algorithm: $R_Z$, $R_Y$, and CNOT \cite{steckmann}. Figure \ref{fig:decomposition_visualization} shows an example of a two-qubit unitary decomposed using these basis gates.
For completeness, we also evaluated each method using its preferred default set of basis gates.
Both Qiskit and Squander use the U3 and CNOT gates as their default basis gate set. The U3 gate is a generic single-qubit rotation gate with 3 Euler angles \cite{qiskit_ugate}.
OpenQL uses $R_Z$, $R_Y$, and CNOT as its default basis gates.
We feel that comparisons are justified even when using different gate sets, as practitioners will often use libraries out of the box without making customizations like changing the set of gates available to the library.

\
\begin{table}[!ht]
    \centering
    \begin{tabular}{lcccccl}\toprule
        \textbf{Method} & \textbf{Qubits} & \textbf{Time (s)} & \textbf{CNOTs} & \textbf{Total Gates} & \textbf{Error} \\ \midrule
        \textbf{Qiskit (Opt = 1)} & 1 & 0.00162 & 0 & 3 & 9E-17 \\
        \textbf{Qiskit (Opt = 2)} & 1 & 0.00235 & 0 & 3 & 9E-17 \\
        \textbf{Qiskit (Opt = 3)} & 1 & 0.00268 & 0 & 3 & 9E-17 \\
        \textbf{OpenQL} & 1 & 0.00122 & 0 & 3 & 6.8E-17 \\
        \textbf{Squander} & 1 & (not supported) & - & - & - \\
        \textbf{Qiskit (Opt = 1)} & 2 & 0.00196 & 3 & 25 & 1.39E-16 \\
        \textbf{Qiskit (Opt = 2)} & 2 & 0.00553 & 3 & 24 & 1.22E-16 \\
        \textbf{Qiskit (Opt = 3)} & 2 & 0.00727 & 3 & 24 & 1.24E-16 \\
        \textbf{OpenQL} & 2 & 0.00140 & 6 & 24 & 1.48E-16 \\
        \textbf{Squander} & 2 & 0.013327 & 2 & 10 & 0.463 \\
        \textbf{Qiskit (Opt = 1)} & 3 & 0.01706 & 20 & 106 & 3.16E-16 \\
        \textbf{Qiskit (Opt = 2)} & 3 & 0.03479 & 20 & 101 & 2.86E-16 \\
        \textbf{Qiskit (Opt = 3)} & 3 & 0.03702 & 20 & 101 & 2.86E-16 \\
        \textbf{OpenQL} & 3 & 0.00300 & 36 & 120 & 2.23E-16 \\
        \textbf{Squander} & 3 & 0.01510 & 40 & 160 & 0.734 \\
        \textbf{Qiskit (Opt = 1)} & 4 & 0.06866 & 100 & 454 & 5.15E-16 \\
        \textbf{Qiskit (Opt = 2)} & 4 & 0.14413 & 100 & 435 & 4.93E-16 \\
        \textbf{Qiskit (Opt = 3)} & 4 & 0.15417 & 100 & 435 & 4.93E-16 \\
        \textbf{OpenQL} & 4 & 0.00957 & 168 & 528 & 3.81E-16 \\
        \textbf{Squander} & 4 & 0.02968 & 60 & 220 & 0.087 \\ \bottomrule
    \end{tabular}
    \caption{Performance comparison of 1, 2, 3, and 4-qubit circuits with custom basis gates, measured by averaging over 100 different unitaries of the appropriate size. The version of Squander we used did not support 1-qubit decompositions.}
    \label{tab:custom-gates}
\end{table}

Table \ref{tab:custom-gates} contains the results of our numerical tests when using the same, custom set of gates ($R_Z$, $R_Y$, and CNOT) for each library.
All results are averaged across 100 random unitaries of the appropriate size for the number of qubits listed.
We observe that Qiskit and OpenQL perform almost identically for single-qubit gates, which are relatively straightforward to decompose.
As we move to multi-qubit operations, more significant differences emerge.
Qiskit's performance remains relatively consistent across optimization levels. This suggests that the higher optimization levels may not provide substantial benefits for these specific decompositions, possibly due to the limited number of qubits being tested.
OpenQL demonstrates competitive performance in terms of error rates, often achieving lower errors than Qiskit and also running several times faster than Qiskit (even with the least aggressive optimization level).
However, OpenQL tends to yield decompositions with a higher number of CNOT gates, which could be a disadvantage in scenarios where two-qubit operations are particularly costly or error-prone.
Squander shows interesting trade-offs: it achieves the lowest CNOT and total gate counts, but at the cost of a significantly higher error rate.
Based on our next results, we surmise that Squander is not optimized for the case of using custom basis gates.

We subsequently use each library's default set of basis gates and re-run the same experiments.
Results are shown in Table \ref{tab:default-gates}.
Since our custom basis gate set used for the results in Table \ref{tab:custom-gates} is the same as OpenQL's default basis gate set, the OpenQL results in Table \ref{tab:default-gates} are identical, but we repeat them for ease of exposition.
We note that for the Qiskit results, running with optimization level 1 yielded circuits with a structure of alternating U3 gates and CNOT gates, which is already an optimal pattern that optimization levels 2 and 3 could not improve upon; therefore, when using Qiskit's default basis gates as in Table \ref{tab:default-gates}, we only report results using optimization level 1.
For larger qubit counts, all methods, including Squander, produce results with errors near machine precision, suggesting that each library is able to decompose the random unitary with an almost perfect approximation.
Squander outperforms OpenQL and Qiskit in terms of CNOT and total gate counts, but it does not compare favorably to OpenQL or Qiskit based on wall clock time.
OpenQL produces circuits with more CNOTs and gates than Qiskit, but OpenQL can run several times faster.
Hence further evaluation of OpenQL and Qiskit is motivated by these results to elucidate their relevant performance; and it is interesting to consider whether the slight advantages in gate counts achieved by Squander could also be realized with improvements to the other libraries, without incurring too great of a computational cost.

\begin{table}[!ht]
    \centering
    \begin{tabular}{lcccccl}\toprule
        \textbf{Method} & \textbf{Qubits} & \textbf{Time (s)} & \textbf{CNOTs} & \textbf{Total Gates} & \textbf{Error} \\ \midrule
        \textbf{Qiskit} & 1 & 0.00182 & 0 & 1 & 8E-17 \\
        \textbf{OpenQL} & 1 & 0.00122 & 0 & 3 & 6.8E-17 \\
        \textbf{Squander} & 1 & (not supported) & - & - & - \\
        \textbf{Qiskit} & 2 & 0.00204 & 3 & 11 & 8E-17 \\
        \textbf{OpenQL} & 2 & 0.00140 & 6 & 24 & 1.48E-16 \\
        \textbf{Squander} & 2 & 2.46151 & 3 & 11 & 7.63E-14 \\
        \textbf{Qiskit} & 3 & 0.01593 & 20 & 57 & 1.34E-16 \\
        \textbf{OpenQL} & 3 & 0.00300 & 36 & 120 & 2.23E-16 \\
        \textbf{Squander} & 3 & 32.46883 & 17 & 54 & 3.23E-12 \\
        \textbf{Qiskit} & 4 & 0.06015 & 100 & 265 & 3.12E-16 \\
        \textbf{OpenQL} & 4 & 0.00957 & 168 & 528 & 3.81E-16 \\
        \textbf{Squander} & 4 & 224.17199 & 77 & 235 & 3.95E-11 \\ \bottomrule
    \end{tabular}
    \caption{Performance comparison of 1, 2, 3, and 4-qubit circuits with default basis gates, measured by averaging over 100 different unitaries of appropriate size. The version of Squander we used did not support 1-qubit decompositions.  All Qiskit optimization levels yielded identical circuits.}
    \label{tab:default-gates}
\end{table}

\section{Conclusions \& Future Work}

In this paper, we demonstrated the utility and performance of several numerical algorithms on the task of automatically synthesizing quantum algorithms from example input-output data.
We considered various aspects of the algorithms' performance, including their scalability, convergence properties, and sensitivity to different step sizes (learning rates).
We also explored factors that may impact the convergence of algorithms on this problem, such as under- and over-constraining, the distance between $Y$ and $X$, the quality of the initial guess, and the conditioning of $X$.
As part of our work, we developed a sequential formulation of the studied problem that, due to its favorable scalability, we recommend as the most feasible approach for use in future research.

We then considered factorization of learned unitaries $U$ into sequences of quantum gates, enabling an end-to-end pipeline from input-output examples to a quantum algorithm that realizes such a transformation.
While Qiskit, OpenQL, and Squander were all successful on this task, we found that---at the cost of more time---Squander was able to produce results with the smallest number of CNOT and total gates, while Qiskit and OpenQL ran several orders of magnitude faster, which could be preferable despite yielding higher gate counts.
Nonetheless, these tests were only on systems with 1--4 qubits, and methods may also exhibit different behavior when using other basis gate sets not considered in our tests.

We envision several potential directions to extend the present work.
From a computing standpoint, it would be useful to explore GPU implementations of the algorithms discussed here and to assess how large of a quantum system one could synthesize in a reasonable amount of time on a GPU.
This is a relevant direction as near-term quantum computers have qubit counts exceeding 100.
There are additional numerical considerations that could also be explored, such as studying the effect of noise on various parts of the system (both in obtaining $U$ and subsequently in factoring it).
It would also be interesting to incorporate penalties to encourage the sparsity of $U$, since---at least for short-depth quantum circuits---the product $U$ of all the gates is often not fully dense.
Finally, given that the delta learning rule algorithm considered in this paper is equivalent to a single-layer neural network, it would be a timely avenue of research to go beyond classical numerical algorithms and investigate deep neural network approaches for the present problem.

\section{Acknowledgments}

Hyde was supported by the U.S.\ Department of Defense (DoD) through the National Defense Science \& Engineering Graduate Fellowship (NDSEG) Program.
The authors acknowledge Anne Kuckertz and Aryan Garg, who helped improve our software implementation.

\newpage

\bibliographystyle{plainnat}
\bibliography{paper}

\begin{thebibliography}{49}
\providecommand{\natexlab}[1]{#1}
\providecommand{\url}[1]{\texttt{#1}}
\expandafter\ifx\csname urlstyle\endcsname\relax
  \providecommand{\doi}[1]{doi: #1}\else
  \providecommand{\doi}{doi: \begingroup \urlstyle{rm}\Url}\fi

\bibitem[{Alan Edelman} and Rao(2005)]{alan_edelman_random_2005}
{Alan Edelman} and N.~Raj Rao.
\newblock Random matrix theory.
\newblock \emph{Acta Numerica}, 14:\penalty0 233--297, May 2005.
\newblock \doi{10.1017/S0962492904000236}.

\bibitem[Barenco et~al.(1995)Barenco, Bennett, Cleve, DiVincenzo, Margolus,
  Shor, Sleator, Smolin, and Weinfurter]{PhysRevA.52.3457}
Adriano Barenco, Charles~H. Bennett, Richard Cleve, David~P. DiVincenzo, Norman
  Margolus, Peter Shor, Tycho Sleator, John~A. Smolin, and Harald Weinfurter.
\newblock Elementary gates for quantum computation.
\newblock \emph{Physical Review A}, 52:\penalty0 3457--3467, Nov 1995.
\newblock \doi{10.1103/PhysRevA.52.3457}.
\newblock URL \url{https://link.aps.org/doi/10.1103/PhysRevA.52.3457}.

\bibitem[Dawson and Nielsen(2005)]{dawson2005solovay}
Christopher~M Dawson and Michael~A Nielsen.
\newblock The {Solovay}-{Kitaev} algorithm, 2005.
\newblock URL \url{https://arxiv.org/abs/quant-ph/0505030}.

\bibitem[Faber et~al.(2003)Faber, Thess, and Giraldi]{giraldi_ga}
Jean Faber, Ricardo~N. Thess, and Gilson~A. Giraldi.
\newblock Learning linear operators by genetic algorithms.
\newblock Technical report, LNCC---National Laboratory for Scientific
  Computing, Rio de Janeiro, Brazil, February 2003.
\newblock URL \url{http://qubit.lncc.br/files/jfaber_Learning-GA.pdf}.

\bibitem[Gilchrist et~al.(2005)Gilchrist, Langford, and
  Nielsen]{Gilchrist_2005}
Alexei Gilchrist, Nathan~K. Langford, and Michael~A. Nielsen.
\newblock Distance measures to compare real and ideal quantum processes.
\newblock \emph{Physical Review A}, 71\penalty0 (6), June 2005.
\newblock ISSN 1094-1622.
\newblock \doi{10.1103/physreva.71.062310}.
\newblock URL \url{http://dx.doi.org/10.1103/PhysRevA.71.062310}.

\bibitem[Giraldi et~al.(2004)Giraldi, Portugal, and Thess]{giraldi2004genetic}
Gilson~A. Giraldi, Renato Portugal, and Ricardo~N. Thess.
\newblock Genetic algorithms and quantum computation, 2004.
\newblock URL \url{https://arxiv.org/abs/cs/0403003}.

\bibitem[{Goubault de Brugière} et~al.(2020){Goubault de Brugière}, Baboulin,
  Valiron, and Allouche]{GOUBAULTDEBRUGIERE2020107001}
Timothée {Goubault de Brugière}, Marc Baboulin, Benoît Valiron, and Cyril
  Allouche.
\newblock Quantum circuits synthesis using householder transformations.
\newblock \emph{Computer Physics Communications}, 248:\penalty0 107001, 2020.
\newblock ISSN 0010-4655.
\newblock \doi{https://doi.org/10.1016/j.cpc.2019.107001}.
\newblock URL
  \url{https://www.sciencedirect.com/science/article/pii/S0010465519303388}.

\bibitem[Greenaway et~al.(2021)Greenaway, Sauvage, Khosla, and
  Mintert]{PhysRevResearch.3.033031}
Sean Greenaway, Fr\'ed\'eric Sauvage, Kiran~E. Khosla, and Florian Mintert.
\newblock Efficient assessment of process fidelity.
\newblock \emph{Physical Review Research}, 3:\penalty0 033031, Jul 2021.
\newblock \doi{10.1103/PhysRevResearch.3.033031}.
\newblock URL \url{https://link.aps.org/doi/10.1103/PhysRevResearch.3.033031}.

\bibitem[Greenberger et~al.(1989)Greenberger, Horne, and
  Zeilinger]{greenberger_going_1989}
Daniel Greenberger, Michael Horne, and Anton Zeilinger.
\newblock Going {Beyond} {Bell}{\textquoteright}s {Theorem}.
\newblock In \emph{Bell{\textquoteright}s {Theorem}, {Quantum} {Theory} and
  {Conceptions} of the {Universe}}, volume~37 of \emph{Fundamental {Theories}
  of {Physics}}, pages 69--72. Springer Netherlands, 1989.
\newblock ISBN 978-90-481-4058-9.

\bibitem[Griffiths(2010)]{qmc062}
Robert~B. Griffiths.
\newblock Unitary time development, September 2010.
\newblock URL
  \url{https://euler.phys.cmu.edu/widom/teaching/33-755/qmc062.pdf}.

\bibitem[{IBM Quantum}(2023{\natexlab{a}})]{qiskit_1q_optimization}
{IBM Quantum}.
\newblock {Optimize1qGates} - {Qiskit}.
\newblock
  \url{https://docs.quantum.ibm.com/api/qiskit/qiskit.transpiler.passes.Optimize1qGates},
  2023{\natexlab{a}}.
\newblock Accessed: 2024-08-26.

\bibitem[{IBM Quantum}(2023{\natexlab{b}})]{qiskit_commutative_cancellation}
{IBM Quantum}.
\newblock {CommutativeCancellation} - {Qiskit}.
\newblock
  \url{https://docs.quantum.ibm.com/api/qiskit/qiskit.transpiler.passes.CommutativeCancellation},
  2023{\natexlab{b}}.
\newblock Accessed: 2024-08-26.

\bibitem[{IBM Quantum}(2023{\natexlab{c}})]{qiskit_inverse_cancellation}
{IBM Quantum}.
\newblock {InverseCancellation} - {Qiskit}.
\newblock
  \url{https://docs.quantum.ibm.com/api/qiskit/qiskit.transpiler.passes.InverseCancellation},
  2023{\natexlab{c}}.
\newblock Accessed: 2024-08-26.

\bibitem[{IBM Quantum}(2023{\natexlab{d}})]{qiskit_ugate}
{IBM Quantum}.
\newblock Ugate - {Qiskit}.
\newblock
  \url{https://docs.quantum.ibm.com/api/qiskit/qiskit.circuit.library.UGate},
  2023{\natexlab{d}}.
\newblock Accessed: 2024-08-23.

\bibitem[Iten et~al.(2016)Iten, Colbeck, Kukuljan, Home, and
  Christandl]{Iten_2016}
Raban Iten, Roger Colbeck, Ivan Kukuljan, Jonathan Home, and Matthias
  Christandl.
\newblock Quantum circuits for isometries.
\newblock \emph{Physical Review A}, 93\penalty0 (3), March 2016.
\newblock ISSN 2469-9934.
\newblock \doi{10.1103/physreva.93.032318}.
\newblock URL \url{http://dx.doi.org/10.1103/PhysRevA.93.032318}.

\bibitem[Javadi-Abhari et~al.(2024)Javadi-Abhari, Treinish, Krsulich, Wood,
  Lishman, Gacon, Martiel, Nation, Bishop, Cross, Johnson, and
  Gambetta]{qiskit2024}
Ali Javadi-Abhari, Matthew Treinish, Kevin Krsulich, Christopher~J. Wood, Jake
  Lishman, Julien Gacon, Simon Martiel, Paul~D. Nation, Lev~S. Bishop,
  Andrew~W. Cross, Blake~R. Johnson, and Jay~M. Gambetta.
\newblock Quantum computing with {Q}iskit, 2024.
\newblock URL \url{https://arxiv.org/abs/2405.08810}.

\bibitem[Khammassi et~al.(2021)Khammassi, Ashraf, Someren, Nane, Krol, Rol,
  Lao, Bertels, and Almudever]{openql2021}
N.~Khammassi, I.~Ashraf, J.~V. Someren, R.~Nane, A.~M. Krol, M.~A. Rol, L.~Lao,
  K.~Bertels, and C.~G. Almudever.
\newblock {OpenQL}: A portable quantum programming framework for quantum
  accelerators.
\newblock \emph{Journal on Emerging Technologies in Computing Systems},
  18\penalty0 (1), Dec. 2021.
\newblock ISSN 1550-4832.
\newblock \doi{10.1145/3474222}.
\newblock URL \url{https://doi.org/10.1145/3474222}.

\bibitem[Kitaev et~al.(2002)Kitaev, Shen, and Vyalyi]{kitaev2002classical}
Alexei~Yu Kitaev, Alexander Shen, and Mikhail~N. Vyalyi.
\newblock \emph{Classical and quantum computation}.
\newblock Number~47 in Graduate Studies in Mathematics. American Mathematical
  Soc., 2002.

\bibitem[Krol et~al.(2022)Krol, Sarkar, Ashraf, Al-Ars, and
  Bertels]{krol2022efficient}
Anna~M. Krol, Aritra Sarkar, Imran Ashraf, Zaid Al-Ars, and Koen Bertels.
\newblock Efficient decomposition of unitary matrices in quantum circuit
  compilers.
\newblock \emph{Applied Sciences}, 12\penalty0 (2):\penalty0 759, 2022.

\bibitem[Liu and Jiang(2014)]{liu2014learning}
Ding Liu and Minghu Jiang.
\newblock Learning quantum operator by quantum adiabatic computation.
\newblock In \emph{2014 12th International Conference on Signal Processing
  (ICSP)}, pages 63--67. IEEE, 2014.

\bibitem[Madden and Simonetto(2022)]{10.1145/3505181}
Liam Madden and Andrea Simonetto.
\newblock Best approximate quantum compiling problems.
\newblock \emph{ACM Transactions on Quantum Computing}, 3\penalty0 (2), Mar.
  2022.
\newblock ISSN 2643-6809.
\newblock \doi{10.1145/3505181}.
\newblock URL \url{https://doi.org/10.1145/3505181}.

\bibitem[Mahmud et~al.(2021)Mahmud, MacGillivray, Chaudhary, and
  El-Araby]{mahmud2021optimizing}
Naveed Mahmud, Andrew MacGillivray, Manu Chaudhary, and Esam El-Araby.
\newblock Optimizing quantum circuits for arbitrary state synthesis and
  initialization.
\newblock In \emph{2021 IEEE 34th International System-on-Chip Conference
  (SOCC)}, pages 19--24. IEEE, 2021.
\newblock \doi{10.1109/SOCC52499.2021.9739614}.

\bibitem[M\"{o}tt\"{o}nen et~al.(2005)M\"{o}tt\"{o}nen, Vartiainen, Bergholm,
  and Salomaa]{10.5555/2011670.2011675}
Mikko M\"{o}tt\"{o}nen, Juha~J. Vartiainen, Ville Bergholm, and Martti~M.
  Salomaa.
\newblock Transformation of quantum states using uniformly controlled
  rotations.
\newblock \emph{Quantum Informatino and Computation}, 5\penalty0 (6):\penalty0
  467–473, Sep. 2005.
\newblock ISSN 1533-7146.

\bibitem[Murnaghan(1958)]{murnaghan1958orthogonal}
Francis~Dominic Murnaghan.
\newblock The orthogonal and symplectic groups.
\newblock \emph{Communications of the Dublin Institute for Advanced Studies},
  1958.

\bibitem[Nielsen and Chuang(2000)]{nielsen_quantum_2000}
Michael Nielsen and Isaac Chuang.
\newblock \emph{Quantum {Computation} and {Quantum} {Information}}.
\newblock Cambridge University Press, Cambridge, UK, 2000.
\newblock ISBN 0-521-63235-8.

\bibitem[Nocedal and Wright(2006)]{wright_numerical_2006}
Jorge Nocedal and Stephen Wright.
\newblock \emph{Numerical {Optimization}}.
\newblock Springer {Series} in {Operations} {Research} and {Financial}
  {Engineering}. Springer-Verlag, New York, second edition, 2006.
\newblock ISBN 978-0-387-30303-1.

\bibitem[Ozols(2009)]{ozols_how_2009}
Maris Ozols.
\newblock How to generate a random unitary matrix, March 2009.
\newblock URL
  \url{http://home.lu.lv/~sd20008/papers/essays/Random%20unitary%20%5Bpaper%5D.pdf}.

\bibitem[Paige and Wei(1994)]{paige1994history}
Christopher~C. Paige and Musheng Wei.
\newblock History and generality of the cs decomposition.
\newblock \emph{Linear Algebra and its Applications}, 208:\penalty0 303--326,
  1994.

\bibitem[Rakyta and Zimbor{\'a}s(2022{\natexlab{a}})]{rakyta2022approaching}
P{\'e}ter Rakyta and Zolt{\'a}n Zimbor{\'a}s.
\newblock Approaching the theoretical limit in quantum gate decomposition.
\newblock \emph{Quantum}, 6:\penalty0 710, 2022{\natexlab{a}}.

\bibitem[Rakyta and Zimbor{\'a}s(2022{\natexlab{b}})]{rakyta2022efficient}
P{\'e}ter Rakyta and Zolt{\'a}n Zimbor{\'a}s.
\newblock Efficient quantum gate decomposition via adaptive circuit
  compression, 2022{\natexlab{b}}.
\newblock URL \url{https://arxiv.org/abs/2203.04426}.

\bibitem[Rakyta et~al.(2022)Rakyta, Morse, N{\'a}dori, Majnay-Tak{\'a}cs,
  Mencer, and Zimbor{\'a}s]{rakyta2022highly}
Peter Rakyta, Gregory Morse, Jakab N{\'a}dori, Zita Majnay-Tak{\'a}cs, Oskar
  Mencer, and Zolt{\'a}n Zimbor{\'a}s.
\newblock Highly optimized quantum circuits synthesized via data-flow engines,
  2022.
\newblock URL \url{https://arxiv.org/abs/2211.07685}.

\bibitem[Reck et~al.(1994)Reck, Zeilinger, Bernstein, and
  Bertani]{reck1994experimental}
Michael Reck, Anton Zeilinger, Herbert~J Bernstein, and Philip Bertani.
\newblock Experimental realization of any discrete unitary operator.
\newblock \emph{Physical Review Letters}, 73\penalty0 (1):\penalty0 58, 1994.

\bibitem[Rudolph et~al.(2022)Rudolph, Chen, Miller, Acharya, and
  Perdomo-Ortiz]{rudolph2022decomposition}
Manuel~S. Rudolph, Jing Chen, Jacob Miller, Atithi Acharya, and Alejandro
  Perdomo-Ortiz.
\newblock Decomposition of matrix product states into shallow quantum circuits,
  2022.
\newblock URL \url{https://arxiv.org/abs/2209.00595}.

\bibitem[Sanderson(2010)]{sanderson_armadillo:_2010}
Conrad Sanderson.
\newblock Armadillo: {An} {Open} {Source} {C}++ {Linear} {Algebra} {Library}
  for {Fast} {Prototyping} and {Computationally} {Intensive} {Experiments}.
\newblock Technical {Report}, NICTA, 2010.

\bibitem[Sanderson and Curtin(2019)]{sanderson2019practical}
Conrad Sanderson and Ryan Curtin.
\newblock Practical sparse matrices in c++ with hybrid storage and
  template-based expression optimisation.
\newblock \emph{Mathematical and Computational Applications}, 24\penalty0
  (3):\penalty0 70, 2019.

\bibitem[Sch{\"o}nemann(1966)]{schonemann_generalized_1966}
Peter Sch{\"o}nemann.
\newblock A generalized solution of the orthogonal procrustes problem.
\newblock \emph{Psychometrika}, 31\penalty0 (1):\penalty0 1--10, March 1966.
\newblock \doi{10.1007/BF02289451}.

\bibitem[Shende et~al.(2005)Shende, Bullock, and Markov]{shende2005synthesis}
Vivek~V. Shende, Stephen~S. Bullock, and Igor~L. Markov.
\newblock Synthesis of quantum logic circuits.
\newblock In \emph{Proceedings of the 2005 Asia and South Pacific Design
  Automation Conference}, pages 272--275, 2005.
\newblock \doi{10.1145/1120725.1120847}.

\bibitem[Shor(1997)]{shor_polynomial-time_1997}
Peter Shor.
\newblock Polynomial-{Time} {Algorithms} for {Prime} {Factorization} and
  {Discrete} {Logarithms} on a {Quantum} {Computer}.
\newblock \emph{SIAM Journal on Computing}, 26\penalty0 (5):\penalty0
  1484--1509, October 1997.
\newblock \doi{10.1137/S0097539795293172}.

\bibitem[Song and Williams(2003)]{song2003computational}
Lin Song and Colin~P. Williams.
\newblock Computational synthesis of any n-qubit pure or mixed state.
\newblock In \emph{Quantum Information and Computation}, volume 5105, pages
  195--203. SPIE, 2003.

\bibitem[Steckmann(2020)]{steckmann}
Thomas Steckmann.
\newblock Unitary decompositions for quantum circuits, 2020.
\newblock URL
  \url{https://thomassteckmann.com/Universal_Gate_Decomposition.pdf}.

\bibitem[Stewart(1982)]{stewart1982computing}
Gilbert~W. Stewart.
\newblock Computing the {CS} decomposition of a partitioned orthonormal matrix.
\newblock \emph{Numerische Mathematik}, 40\penalty0 (3):\penalty0 297--306,
  1982.

\bibitem[Toronto and Ventura(2006)]{toronto_learning_2006}
Neil Toronto and Dan Ventura.
\newblock Learning {Quantum} {Operators} {From} {Quantum} {State} {Pairs}.
\newblock In \emph{2006 {IEEE} International Conference on Evolutionary
  Computation}, pages 2607--2612, Vancouver, BC, Canada, 2006. IEEE.
\newblock \doi{10.1109/CEC.2006.1688634}.

\bibitem[Tucci(1998)]{tucci1998rudimentary}
Robert~R. Tucci.
\newblock A rudimentary quantum compiler, 1998.
\newblock URL \url{https://arxiv.org/abs/quant-ph/9805015}.

\bibitem[Tucci(1999)]{tucci1999rudimentary}
Robert~R. Tucci.
\newblock A rudimentary quantum compiler (2nd ed.), 1999.
\newblock URL \url{https://arxiv.org/abs/quant-ph/9902062}.

\bibitem[Tucci(2005)]{tucci2005KAKDecomp}
Robert~R. Tucci.
\newblock An introduction to {Cartan}'s {KAK} decomposition for {QC}
  programmers, 2005.
\newblock URL \url{https://arxiv.org/abs/quant-ph/0507171}.

\bibitem[Van~Loan(1985)]{van1985computing}
Charles Van~Loan.
\newblock Computing the {CS} and the generalized singular value decompositions.
\newblock \emph{Numerische Mathematik}, 46\penalty0 (4):\penalty0 479--491,
  1985.

\bibitem[Ventura(2000)]{ventura_learning_2000}
Dan Ventura.
\newblock Learning {Quantum} {Operators}.
\newblock \emph{Proceedings of the Joint Conference on Information Sciences},
  pages 750--752, 2000.

\bibitem[Vlachogiannis(2010)]{vlachogiannis2010learning}
John~G. Vlachogiannis.
\newblock Learning linear operators by coordinated aggregation-based {PSO}.
\newblock \emph{Journal of Experimental \& Theoretical Artificial
  Intelligence}, 22\penalty0 (4):\penalty0 311--319, 2010.

\bibitem[Vlachogiannis(2013)]{vlachogiannis2013improved}
John~G. Vlachogiannis.
\newblock An improved linear learning algorithm based on swarm intelligence.
\newblock \emph{Quantum Matter}, 2\penalty0 (3):\penalty0 199--204, 2013.

\end{thebibliography}

\newpage

\appendix

\section{Derivation of the Hessian of the Frobenius Norm Error}\label{sec:Hessian}

Following the results of Toronto and Ventura \citep{toronto_learning_2006}, we have
\begin{align}
(\nabla_U f(U))_{jk} &= e_j^T (U X - Y) X^H e_k,
\intertext{where $e_j$ and $e_k$ are the Cartesian basis vectors,}
&= \Tr(e_j^T (UX - Y) X^H e_k) \\
&= \Tr( (UX-Y)X^H e_k e_j^T),
\intertext{where we have used the cyclic property of the trace,}
&= \Tr( U X X^H e_k e_j^T) - \Tr( Y X^H e_k e_j^T).
\end{align}
The second term is constant, and hence will have zero gradient. Thus,
\begin{equation}
\nabla_U \left[(\nabla_U f(U))_{jk} \right]= \nabla_U \Tr(U X X^H e_k e_j^T).
\end{equation}
If we separate $U = U_R + i U_I$, then, for an arbitrary analytic function $F(U)$, we know
\begin{equation}
\nabla_U F(U) = \frac{1}{2} \left( \nabla_{U_R} F(U) - i \nabla_{U_I} F(U) \right).
\end{equation}
Hence,
\begin{equation}
\begin{aligned}
\nabla_U \left[(\nabla_U f(U))_{jk}\right] = &\frac{1}{2} \nabla_{U_R} \Tr(U X X^H e_k e_j^T)\\
&-i \frac{1}{2} \nabla_{U_I} \Tr(U XX^H e_k e_j^T) \\
= &(XX^H e_k e_j^T)^T \\
= &e_j e_k^T (X^H)^T X^T.
\end{aligned}
\end{equation}
In indicial notation, we have
\begin{equation}\label{eq:OrthogonalProcrustesHessian}
\begin{aligned}
\frac{\partial^2 f(U)}{\partial U_{jk} \partial U_{mn}} &= e_m^T e_j e_k^T (X^H)^T X^T e_n \\
&= \delta_{mj} \delta_{kb} \bar{X}_{bc} X_{dc} \delta_{dn} \\
&= \delta_{jm} \bar{X}_{kc} X_{nc}.
\end{aligned}
\end{equation}
This fourth-order tensor may be rearranged into a matrix for use in Newton's method.

\end{document}